\documentclass[10pt]{article}
\usepackage{amsmath}
\usepackage{color}
\usepackage{graphicx}
\usepackage{fancyhdr}
\usepackage[english]{babel}
\usepackage[autostyle]{csquotes}
\usepackage{verbatim}
\usepackage{amsfonts,longtable}
\usepackage{epsfig}
\usepackage{mathrsfs}
\usepackage{mathtools}
\usepackage{amsfonts}
\usepackage{color}

\def\R{\hbox{{\rm I}\kern-0.2em{\rm R}\kern0.2em}}

\def\d{{\rm d}}

\def\bn{\begin{equation}}
\def\en{\end{equation}}
\def\bny{\begin{eqnarray}}
\def\eny{\end{eqnarray}}
\def\be{\begin{eqnarray*}}
\def\ee{\end{eqnarray*}}
\def\bc{\begin{center}}
\def\ec{\end{center}}
 
\def\p{\partial}
\def\({\left(}
\def\){\right  )}
\def\[{\left[}
\def\]{\right]}
\def\bc{\begin{center}}
\def\ec{\end{center}}

\newtheorem{dfn}{Definition}[section]
\newtheorem{thm}{Theorem}[section]
\newtheorem{rem}{Remark}[section]
\newtheorem{pro}{Proposition}[section]

\newtheorem{cor}{Corollary}[section]
\newtheorem{lem}{Lemma}[section]
\newtheorem{exm}{Example}[section]

\def\bn{\begin{equation}}
\def\en{\end{equation}}
\def\bny{\begin{eqnarray}}
\def\eny{\end{eqnarray}}
\def\ba{\begin{array}{lllll}}
\def\ea{\end{array}}
\def\be{\begin{eqnarray*}}
\def\ee{\end{eqnarray*}}
\def\bdn{\begin{dfn}}
\def\edn{\end{dfn}}
\def\btm{\begin{thm}}
\def\etm{\end{thm}}
\def\bpf{\begin{proof}}
\def\epf{\end{proof}}
\def\bpn{\begin{pro}}
\def\epn{\end{pro}}
\def\brk{\begin{rem}}
\def\erk{\end{rem}}
\def\bcy{\begin{cor}}
\def\ecy{\end{cor}}
\def\blm{\begin{lem}}\def\elm{\end{lem}}
\def\bex{\begin{exm}}
\def\eex{\end{exm}}

 \def\R{\mathscr{R}}
\def\Q{{\cal Q}}

\begin{document}
\bc \Large{{\bf
Symmetry, reductions and exact solutions  of  systems of difference equations.
}}\ec
\medskip

\bc
JJ Bashingwa, AH Kara and M Folly-Gbetoula\\ \vspace{1cm}
School of Mathematics, University of the Witwatersrand, \\
Johannesburg, Wits 2001, South Africa.\ec
%
%
%
%

\begin{abstract}
\noindent We apply symmetry and invariance methods to analyse systems of difference equations. Non trivial symmetries are derived and their exact solutions obtained.
\end{abstract}
\section{Introduction}
The use of symmetry methods for solving difference equations ($\Delta$Es) have been growing rapidly in recent years \cite{maeda,hyd,levi,qui,mensah}.
Maeda \cite{maeda} showed how to use symmetry methods to simplify and obtained    solutions to autonomous system of first-order ordinary difference equations (O$\Delta$Es). In \cite{levi,qui}, the authors extended Maeda's ideas and presented a series-based methods for obtaining the linearized symmetry condition (LSC).
In \cite{hyd}, the author introduced a method for obtaining symmetries (in closed form) and first integrals of
O$\Delta$Es. \par \noindent
In this paper, we extend these approaches to  system of difference equations (S$\Delta$Es). We generalize some results which have been found by some authors.
Elsayed \cite{elsa} studied the solutions of the following system
$$x_{n+1}=\frac{x_{n-1}}{\pm 1+ y_n x_{n-1}}, \quad y_{n+1}=\frac{y_{n-1}}{\mp 1+ x_n y_{n-1}},$$
where $y_0,\,y_{-1},\,x_0,\,x_{-1} \neq 0.$
Kurbanli et al. \cite{kur} investigated the positive solutions of
$$x_{n+1}=\frac{x_{n-1}}{ 1+ y_n x_{n-1}}, \quad y_{n+1}=\frac{y_{n-1}}{ 1+ x_n y_{n-1}},$$
where $y_0,\,y_{-1},\,x_0,\,x_{-1} \in [0, \infty )$.
Elsayed and Ibrahim \cite{elsaibra} studied the solutions of the system
$$x_{n+1} = \frac{x_{n-2}y_{n-1}}{y_n(\pm 1\pm  x_{n-2} y_{n-1})} ,\quad y_{n+1} = \frac{y_{n-2}x_{n-1}}{x_n(\pm1 \pm y_{n-2} x_{n-1})},$$
where initial conditions  $x_{-2} ,\, x_{-1},\,x_0 ,\,y_{-2},\,y_{-1}$ and $y_0$  are nonzero and real.\par \noindent
In our first application, we generalize the results in \cite{elsa} and \cite{kur}. The system under consideration is
\bn x_{n+1} = \frac{x_{n-1}}{a + x_{n-1} y_{n}} ,\quad y_{n+1} = \frac{y_{n-1}}{b + y_{n-1} x_{n}},\label{eq1} \en
where $a$ and $b$ are real numbers. In the second application, we generalize the results in \cite{elsaibra} and study the system
\bn x_{n+1} = \frac{x_{n-2}y_{n-1}}{y_n(a + b x_{n-2} y_{n-1})} ,\quad y_{n+1} = \frac{y_{n-2}x_{n-1}}{x_n(c + d y_{n-2} x_{n-1})},\label{eq111} \en
where $a,b,c$ and $d$ are real numbers. \par \noindent
In order to use Lie symmetry method we \lq shift \rq both
 \eqref{eq1} and  \eqref{eq111}. We  study the systems
 \bn u_{n+2} = \frac{u_n}{a + u_n v_{n+1}} ,\quad v_{n+2} = \frac{v_n}{b + v_n u_{n+1}} \label{eq2}\en
and
\bn x_{n+3} = \frac{x_ny_{n+1}}{y_{n+2}(a + bx_n y_{n+1})} ,\quad y_{n+3} = \frac{y_nx_{n+1}}{x_{n+2}(c + dy_n x_{n+1})} \label{eq112}\en instead.
\section{Overview of Lie analysis of systems of difference equations}
Let \cite{saha1, LW1}
\bn
 x_{n+N}^i=\omega_i(n,x_n ^1,...,x_n^r,...,
 x_{n+N-1}^1,...,x_{n+N-1}^r),\quad i=1,...,r,\label{eq3}\en
where $\omega_i$ are   such that $\frac{\p \omega_i	}{\p x_n^i}\neq 0 $,
be an $N$-th   order system of $r$ difference equations (S$\Delta$Es). Consider the transformation
$$\Gamma : (n,x_n ^1,...,x_n^r,...,
 x_{n+N-1}^1,...,x_{n+N-1}^r)\mapsto (\bar{n},\bar{x}_n ^1,...,\bar{x}_n^r,...,
 \bar{x}_{n+N-1}^1,...,\bar{x}_{n+N-1}^r).$$ We seek for a one-parameter (local) Lie group of transformations
\bn  \begin{array}{lll}
 \bar{x}_{n+j}^i=x_{n+j}^i + \epsilon  S^jQ_i(n,x_n ^1,...,x_n^r,...,
 x_{n+N-1}^1,...,x_{n+N-1}^r)+O(\epsilon^2)\end{array}, \en
where  $Q_i$, $i=1,...,r$, $j=1,...,N$, are continuous functions which we shall refer to as  characteristics. $S$ is the "shift" operator and it is defined as follows
\bn  S: n \mapsto n+1 ,\quad
S^k (x_n^i) =x_{n+k}^i .  \en
The symmetry condition for the S$\Delta$Es \eqref{eq3} is
\bn
 \bar{x}_{n+N}^i=\omega_i(\bar{n},\bar{x}_n ^1,...,\bar{x}_n^r,\bar{x}_{n+1}^1,...,
 \bar{x}_{n+1}^r,...,
 \bar{x}_{n+N-1}^1,...,\bar{x}_{n+N-1}^r), \quad i=1,...,r,\label{eq66}\en
whenever \eqref{eq3} holds. Lie symmetries are obtained by linearizing the symmetry condition \eqref{eq66} about the identity. We have the following system of linearized symmetry condition (SLSC)
\bn
  S^N Q_i -X\omega_i=0, \quad i=1,...,r, \label{eq4}\en
where the symmetry generator $X$ is given by
\bn X= \sum\limits_{j=0}^{N-1}\left(\sum\limits_{i=1}^{r}S^jQ_i
 \frac{\p}{\p x_{n+j}^i}\right).\en
\bdn A function $w_n$  is invariant function under the Lie group of transformations $\Gamma$ if \edn
\bn X(w_n)=0,\en
where $w_n$ can be found by solving the characteristic equation
\bn\frac{\d x_n^1}{Q_1} =...=\frac{\d x_{n}^r}{Q_r}=....=\frac{\d x_{n+N-1}^1}{S^{N-1}Q_1}=...=\frac{\d x_{n+N-1}^r}{S^{N-1}Q_r}. \en
{\bf Note:} We are interested on {\it point symmetry}, the characteristics are of the form $Q_i=Q_i(n,x_n^1,...,x_n^r)$.\par \noindent
Substituting \eqref{eq3} in SLSC \eqref{eq4} we obtain a system of difference equations containing  functions, $Q_n, \quad n=1,...,N$  with different arguments. This characteristic functions can be solved after a series of steps:\\
 $\bullet$  We check if  { the Jacobian } of $\omega$ is not zero, i.e.,
\bn \left|\frac{\p (\omega_1,...,\omega_r)}{(x_n^1,...,x_n^r)} \right|=\left|\ba \frac{\p \omega_1}{x_n^1}&...&\frac{\p \omega_1}{x_n^r}\\\vdots &&\vdots\\\frac{\p \omega_r}{x_n^1}&...&\frac{\p \omega_r}{x_n^r} \ea \right|\neq 0.\en This condition is always satisfied as we have $\frac{\omega_i}{x_n^i}\neq 0$.  Then, we can apply the implicit function theorem \cite{kh} by writing
$$x_n^i=\psi_i(n,x_{n+1}^1,...,
 x_{n+1}^r,...,
 x_{n+N-1}^1,...,x_{n+N-1}^r,\omega_i),\quad i=1,...,r,$$ where $\psi_i$ are differentiable functions. Furthermore, for $k=1,...,r$ and $j=1,...,N-1$,  we have

 $$\frac{\p \psi_i}{\p x_{n+j}^k}=-\frac{\frac{\omega_i(n,x_n ^1,...,x_n^r,x_{n+1}^1,...,
 x_{n+1}^r,...,
 x_{n+N-1}^1,...,x_{n+N-1}^r)}{x_{n+j}^k}}{\frac{\omega_i(n,x_n ^1,...,x_n^r,x_{n+1}^1,...,
 x_{n+1}^r,...,
 x_{n+N-1}^1,...,x_{n+N-1}^r)}{x_n^i}}.$$
$\bullet$ We differentiate the system as follows: the $i$-th equation is differentiated with respect to $x_n^i$ keeping $\omega_i$ fixed.
In this way we eliminate in the system the functions    $Q_i(n+N,\omega_1,...,\omega_r),\quad i=1,...,r$. If there is need to eliminate functions $Q_i(n+N-1,x_{n+N-1}^1,...,x_{n+N-1}^r),\quad i=1,...,r,$ in the system, we multiply each equation by whatever factor is needed to make the coefficient of $Q_i(n+N-1,x_{n+N-1}^1,...,x_{n+N-1}^r),\quad i=1,...,r, $ to be 1, then differentiate the system as before, i.e, differentiate the $i$-th equation with respect to $x_n^i$ (keeping $x_{n+N-1}^i$ fixed). The process is repeated until all undesirable terms are removed. Finally, the resulting partially differential equations PDEs  
 are split into a system of PDEs by grouping together all terms with the same dependency on
$x_{n+1}^1,...,x_{n+1}^r,...,x_{n+N-1}^1,...,x_{n+N-1}^r$.
The resulting system can be solved for the $Q_i$'s.
The last step consist of substituting the characteristics found in the SLSC to obtain their explicit form.
\section{Example 1}
Consider the second order system of difference equations \eqref{eq2} \cite{elsa, kur}.  Imposing the SLSC \eqref{eq4}, we get
\bn \begin{array}{ll}
Q_1(n+2,\omega_1,\omega_2)- \frac{a Q_1(n,u_n,v_n)}{(a+u_nv_{n+1})^2}+\frac{u_n^2 Q_2(n+1,u_{n+1},v_{n+1})}{(a+u_nv_{n+1})^2}=0,\\Q_2(n+2,\omega_1,\omega_2)- \frac{b Q_2(n,u_n,v_n)}{(b+v_nu_{n+1})^2}+\frac{v_n^2 Q_1(n+1,u_{n+1},v_{n+1})}{(b+v_nu_{n+1})^2}=0.
\end{array}   \label{eq5}\en
To simplify our work, we let
$Q_1=Q_1(n,u_n)$ and $Q_2=Q_2(n,v_n).$
The system \eqref{eq5} becomes
\begin{eqnarray}
Q_1(n+2,\omega_1)- \frac{a Q_1(n,u_n)}{(a+u_nv_{n+1})^2}+\frac{u_n^2 Q_2(n+1,v_{n+1})}{(a+u_nv_{n+1})^2}=0,\label{eq6}\\Q_2(n+2,\omega_2)- \frac{b Q_2(n,v_n)}{(b+v_nu_{n+1})^2}+\frac{v_n^2 Q_1(n+1,u_{n+1})}{(b+v_nu_{n+1})^2}=0.
   \label{eq7}\end{eqnarray}
We then proceed by differentiating \eqref{eq6} and \eqref{eq7} with respect to $u_n$ and $v_n$ respectively  (keeping $\omega_1$ and $\omega_2$ fixed). This leads to
\begin{eqnarray}
-a Q_1'(n,u_n)+a Q_2'(n+1,v_{n+1})+ \frac{2a Q_1(n,u_n)}{u_n}=0,\label{eq8}\\-b Q_2'(n,v_n)+b Q_1'(n+1,u_{n+1})+ \frac{2b Q_2(n,v_n)}{v_n}=0,\label{eq9}
\end{eqnarray}
where the prime denotes derivative with respect to the continuous variable. Differentiating \eqref{eq8} and \eqref{eq9} with respect to $u_n$ and $v_n$ respectively  (keeping $u_{n+1}$ and $v_{n+1}$ fixed), we obtain
\begin{eqnarray}
-a Q_1''(n,u_n)-\frac{2a Q_1(n,u_n)}{u_n^2}+\frac{2a Q_1'(n,u_n)}{u_n}=0,\\-b Q_2''(n,v_n)-\frac{2b Q_2(n,v_n)}{v_n^2}+\frac{2b Q_2'(n,v_n)}{v_n}=0
\end{eqnarray}
whose solutions are given by
\begin{eqnarray}
Q_1(n,u_n)=F_1(n)u_n+ F_2(n)u_n^2 \;\text{ and } \;
Q_2(n,u_n)=F_3(n)v_n+ F_4(n)v_n^2\label{eq11}.
\end{eqnarray}
The last step will consist of substituting \eqref{eq11} in \eqref{eq8} \& \eqref{eq9} to get the characteristics given by
\begin{eqnarray}
Q_1=(C_2(-1)^n-C_1 )u_n \text{ and } Q_2=(C_1 +C_2(-1)^n)v_n.
\end{eqnarray}
The generators of the Lie point symmetry are
\bn
\begin{array}{ll}
X_1=-u_n \p_{u_n}+v_n \p_{v_n},\quad
X_2=(-1)^nu_n \p_{u_n}+(-1)^nv_n \p_{v_n}.
\end{array}\label{eq21}
\en
We can   easily   check that the functions
\bny w_n=v_n u_{n+1}\label{eq22}\\z_n=u_nv_{n+1} \label{eq23}\eny
are invariants under $X_2$ given in \eqref{eq21} and that
\bn \begin{array}{ll}
 w_{n+1}=\frac{z_n}{a+z_n},\quad z_{n+1}=\frac{w_n}{b +w_n} . \end{array}\label{eq24}\en
Let \bn \ba z_n=\frac{1}{T_n} \quad  \text { and } \quad w_n=\frac{1}{S_n}.\label{tras}\ea\en
%
Using the transformations \eqref{tras}, the  system
\eqref{eq24} becomes linear
\bn \ba S_{n+1}=a T_n +1,\quad T_{n+1}=b S_n+1 .\ea \label{eq25} \en
The general solution of \eqref{eq25} is given by
\bn  \ba S_n=&\left[  \frac{a^{\frac{n}{2}}b^{\frac{n}{2}}+(-1)^na^{\frac{n}{2}}b^{\frac{n}{2}}}{2}\right] S_0 +\left[  \frac{a^{\frac{n+1}{2}}b^{\frac{n-1}{2}}-(-1)^na^{\frac{n+1}{2}}b^{\frac{n-1}{2}}}{2}\right]T_0 \\& +\frac{\left(  \sum\limits_{i=0}^{\frac{n-2}{2}}a^ib^i+\sum\limits_{i=0}^{\frac{n-2}{2}}a^{i+1}b^i\right)+(-1)^n\left(  \sum\limits_{i=0}^{\frac{n-2}{2}}a^ib^i+\sum\limits_{i=0}^{\frac{n-2}{2}}a^{i+1}b^i\right) }{2}\\&+\frac{\left(  \sum\limits_{i=0}^{\frac{n-1}{2}}a^ib^i+\sum\limits_{i=0}^{\frac{n-3}{2}}a^{i+1}b^i\right)-(-1)^n\left(  \sum\limits_{i=0}^{\frac{n-1}{2}}a^ib^i+\sum\limits_{i=0}^{\frac{n-3}{2}}a^{i+1}b^i\right) }{2},\\
T_n=&\left[  \frac{a^{\frac{n}{2}}b^{\frac{n}{2}}+(-1)^na^{\frac{n}{2}}b^{\frac{n}{2}}}{2}\right] T_0 +\left[  \frac{a^{\frac{n-1}{2}}b^{\frac{n+1}{2}}-(-1)^na^{\frac{n-1}{2}}b^{\frac{n+1}{2}}}{2}\right]S_0 \\
 &+\frac{\left(  \sum\limits_{i=0}^{\frac{n-2}{2}}a^ib^i+\sum\limits_{i=0}^{\frac{n-2}{2}}a^{i}b^{i+1}\right)+(-1)^n\left(  \sum\limits_{i=0}^{\frac{n-2}{2}}a^ib^i+\sum\limits_{i=0}^{\frac{n-2}{2}}a^{i}b^{i+1}\right) }{2}\\&+\frac{\left(  \sum\limits_{i=0}^{\frac{n-1}{2}}a^ib^i+\sum\limits_{i=0}^{\frac{n-3}{2}}a^{i}b^{i+1}\right)-(-1)^n\left(  \sum\limits_{i=0}^{\frac{n-1}{2}}a^ib^i+\sum\limits_{i=0}^{\frac{n-3}{2}}a^{i}b^{i+1}\right) }{2}. \ea \label{eq26}\en
The solutions \eqref{eq26} can be split into
\begin{align} &S_{2n}=a^n b^n S_0 + \sum\limits_{i=0}^{n-1} a^ib^i + \sum\limits_{i=0}^{n-1}a^{i+1}b^i,\;
T_{2n}=a^n b^n T_0+ \sum\limits_{i=0}^{n-1} a^ib^i + \sum\limits_{i=0}^{n-1}a^{i}b^{i+1},\nonumber\\
& S_{2n+1}=a^{n+1} b^n T_0 + \sum\limits_{i=0}^{n} a^ib^i + \sum\limits_{i=0}^{n-1}a^{i+1}b^i,\nonumber\\&
T_{2n+1}=a^n b^{n+1} S_0+ \sum\limits_{i=0}^{n} a^ib^i + \sum\limits_{i=0}^{n-1}a^{i}b^{i+1}.  \label{eq29} \end{align}
Equations \eqref{eq22} and \eqref{eq23} can be written
\bn\ba u_{n+1}=\frac{w_n}{v_n},\quad v_{n+1}=\frac{z_n}{u_n}. \label{eq30} \ea\en
Invoking \eqref{tras}, system \eqref{eq30} becomes
\bn \ba u_{n+1}=\frac{1}{S_n v_n},\quad v_{n+1}=\frac{1}{T_nu_n}, \ea \label{eq31}\en
where the functions $S_n$ and $T_n$ are given in \eqref{eq26}.\\
{\bf Note}: The order of the system \eqref{eq2} has been reduced by one in \eqref{eq31}.\\
The general solution of \eqref{eq31} which is also the general solution of the original system \eqref{eq2} is
\bn \ba u_n=\frac{u_0}{2}\left[ \prod\limits_{r=0}^{\frac{n-2}{2}}\frac{T_{2r}}{S_{2r+1}}+(-1)^n\prod\limits_{r=0}^{\frac{n-2}{2}}\frac{T_{2r}}{S_{2r+1}} \right] +\frac{1}{2v_0}\left[\frac{\prod\limits_{r=0}^{\frac{n-3}{2}}T_{2r+1}}{\prod\limits_{r=0}^{\frac{n-1}{2}}S_{2r}}-(-1)^n \frac{\prod\limits_{r=0}^{\frac{n-3}{2}}T_{2r+1}}{\prod\limits_{r=0}^{\frac{n-1}{2}}S_{2r}}
 \right] \\
 v_n=\frac{v_0}{2}\left[ \prod\limits_{r=0}^{\frac{n-2}{2}}\frac{S_{2r}}{T_{2r+1}}+(-1)^n\prod\limits_{r=0}^{\frac{n-2}{2}}\frac{S_{2r}}{T_{2r+1}} \right] +\frac{1}{2u_0}\left[\frac{\prod\limits_{r=0}^{\frac{n-3}{2}}S_{2r+1}}{\prod\limits_{r=0}^{\frac{n-1}{2}}T_{2r}}-(-1)^n \frac{\prod\limits_{r=0}^{\frac{n-3}{2}}S_{2r+1}}{\prod\limits_{r=0}^{\frac{n-1}{2}}T_{2r}}
 \right] \label{eq32} \ea \en
which can be split into
\bn \ba
u_{2n}=u_0 \prod\limits_{r=0}^{n-1}\frac{T_{2r}}{S_{2r+1}},&v_{2n}=v_0 \prod\limits_{r=0}^{n-1}\frac{S_{2r}}{T_{2r+1}},\\u_{2n+1}=\frac{\prod\limits_{r=0}^{n-1}T_{2r+1}}{v_0\prod\limits_{r=0}^{n}S_{2r}},&
v_{2n+1}=\frac{\prod\limits_{r=0}^{n-1}S_{2r+1}}{u_0\prod\limits_{r=0}^{n}T_{2r}} \label{eq33}\ea \en
 with
 \begin{align} & T_0 \neq \frac{-1}{a^{r+1}b^r} \left[ \sum\limits_{i=0}^{r}a^ib^i + \sum\limits_{i=0}^{r-1}a^{i+1}b^i\right],
 S_0 \neq \frac{-1}{a^{r}b^{r+1}} \left[ \sum\limits_{i=0}^{r}a^ib^i + \sum\limits_{i=0}^{r-1}a^{i}b^{i+1}\right],\nonumber\\ &
 S_0 \neq \frac{-1}{a^rb^r} \left( \sum\limits_{i=0}^{r-1}a^ib^i + \sum\limits_{i=0}^{r-1}a^{i+1}b^i\right) , T_0 \neq \frac{-1}{a^rb^r} \left[ \sum\limits_{i=0}^{r-1}a^ib^i + \sum\limits_{i=0}^{r-1}a^{i}b^{i+1}\right],
  \label{eq34_1} \end{align}
$r\leq n$. The initial conditions
\bn \ba T_0=\frac{1}{z_0}=\frac{1}{u_0v_1} \quad \text{and} \quad S_0=\frac{1}{w_0}=\frac{1}{v_0u_1} \label{eq28_1}\ea\en
are obtained from \eqref{tras}. The restrictions on the system \eqref{eq2} are obtained by bringing together \eqref{eq34_1} and \eqref{eq28_1}.
\subsection{The case $ab\neq 1$}
From  \eqref{eq29}, using the fact that $s_0=\frac{1}{v_0u_1}$ and $t_0=\frac{1}{u_0v_1}$, we obtain
\begin{align} T_{2r}=\frac{a^rb^r [1-ab-u_0v_1(1+b)]+u_0v_1(1+b)}{u_0v_1(1-ab)},\\
S_{2r}=\frac{a^rb^r [1-ab-v_0u_1(1+a)]+v_0u_1(1+a)}{v_0u_1(1-ab)} \\
T_{2r+1}=\frac{a^{r}b^{r+1} [1-ab-v_0u_1(1+a)]+v_0u_1(1+b)}{v_0u_1(1-ab)} ,\\
S_{2r+1}=\frac{a^{r+1}b^r [1-ab-u_0v_1(1+b)]+u_0v_1(1+a)}{u_0v_1(1-ab)}.
 \label{eq34}\end{align}
The solutions of the system \eqref{eq2}, given in \eqref{eq33}, become
\bn\ba u_{2n}=&u_0 \prod\limits_{r=0}^{n-1}\left\lbrace \frac{a^rb^r [1-ab-u_0v_1(1+b)]+u_0v_1(1+b)}{a^{r+1}b^r [1-ab-u_0v_1(1+b)]+u_0v_1(1+a)}\right\rbrace ,\\
v_{2n}=&v_0 \prod\limits_{r=0}^{n-1}\left\lbrace \frac{a^rb^r [1-ab-v_0u_1(1+a)]+v_0u_1(1+a)}{a^{r}b^{r+1} [1-ab-v_0u_1(1+a)]+v_0u_1(1+b)} \right\rbrace , \\
u_{2n+1}=&u_1(1-ab)\frac{\prod\limits_{r=0}^{n-1}\left\lbrace a^{r}b^{r+1} [1-ab-v_0u_1(1+a)]+v_0u_1(1+b)\right\rbrace }{\prod\limits_{r=0}^{n}\left\lbrace a^rb^r [1-ab-v_0u_1(1+a)]+v_0u_1(1+a)\right\rbrace },\\
v_{2n+1}=&v_1(1-ab)\frac{\prod\limits_{r=0}^{n-1}\left\lbrace a^{r+1}b^r [1-ab-u_0v_1(1+b)]+u_0v_1(1+a)\right\rbrace }{\prod\limits_{r=0}^{n}\left\lbrace a^rb^r [1-ab-u_0v_1(1+b)]+u_0v_1(1+b)\right\rbrace }.\ea \label{eq35}\en

\begin{itemize}
\item If  $ a=1,b\neq 1$, then the solutions \eqref{eq35} become
\bn\ba u_{2n}=&u_0 \prod\limits_{r=0}^{n-1}\left\lbrace \frac{b^r [1-b-u_0v_1(1+b)]+u_0v_1(1+b)}{b^r [1-b-u_0v_1(1+b)]+2u_0v_1}\right\rbrace ,\\v_{2n}=&v_0 \prod\limits_{r=0}^{n-1}\left\lbrace \frac{b^r [1-b-2v_0u_1]+2v_0u_1}{b^{r+1} [1-b-2v_0u_1]+v_0u_1(1+b)} \right\rbrace , \\
u_{2n+1}=&u_1(1-b)\frac{\prod\limits_{r=0}^{n-1}\left\lbrace b^{r+1} [1-b-2v_0u_1]+v_0u_1(1+b)\right\rbrace }{\prod\limits_{r=0}^{n}\left\lbrace b^r [1-b-2v_0u_1]+2v_0u_1\right\rbrace },\\
v_{2n+1}=&v_1(1-b)\frac{\prod\limits_{r=0}^{n-1}\left\lbrace b^r [1-b-u_0v_1(1+b)]+2u_0v_1\right\rbrace }{\prod\limits_{r=0}^{n}\left\lbrace b^r [1-b-u_0v_1(1+b)]+u_0v_1(1+b)\right\rbrace }.\ea \label{eq36}\en
 For a particular value  $b=-1$,  we obtain
\bn\ba u_{2n}=u_0 \prod\limits_{r=0}^{n-1}\left\lbrace \frac{(-1)^r}{(-1)^r+u_0v_1}\right\rbrace ,&
v_{2n}=v_0 \prod\limits_{r=0}^{n-1}\left\lbrace \frac{(-1)^r [1-v_0u_1]+v_0u_1}{(-1)^{r+1} [1-v_0u_1]} \right\rbrace , \\
u_{2n+1}=u_1\frac{\prod\limits_{r=0}^{n-1}\left\lbrace (-1)^{r+1} [1-v_0u_1]\right\rbrace }{\prod\limits_{r=0}^{n}\left\lbrace (-1)^r [1-v_0u_1]+v_0u_1\right\rbrace },&
v_{2n+1}=v_1\frac{\prod\limits_{r=0}^{n-1}\left\lbrace (-1)^r +u_0v_1\right\rbrace }{\prod\limits_{r=0}^{n}\left\lbrace (-1)^r \right\rbrace }.\ea \label{eq37}\en
The solutions in \eqref{eq37} can be split as follows
\bn \ba u_{4n}&=\frac{u_0}{(1-u_0^2v_1^2)^n},\qquad  &u_{4n+1}&=\frac{u_1(1-v_0u_1)^{2n}}{(1-2v_0u_1)^n},\\
u_{4n+2}&=\frac{u_0}{(1+u_0v_1)(1-u_0^2v_1^2)^n}, \qquad
&u_{4n+3}&=\frac{(-1)^n(v_0u_1-1)^{2n+1}}{(2v_0u_1-1)^{n+1}},\\
v_{4n}&=\frac{v_0(1-2v_0u_1)^n}{(1-v_0u_1)^{2n}},\qquad
&v_{4n+1}&=v_1(1-u_0^2v_1^2)^n,\\
v_{4n+2}&=\frac{-v_0(1-2v_0u_1)^n}{(1-v_0u_1)^{2n+1}},\qquad
&v_{4n+3}&=-v_1(1+u_0v_1)(1-u_0^2v_1^2)^n.\label{eq38}\ea\en
If we let $v_{4n}=y_{4n-1},v_{4n+1}=y_{4n},v_{4n+2}=y_{4n+1},v_{4n+3
}=y_{4n+2},u_{4n}=x_{4n-1},u_{4n+1}=x_{4n},u_{4n+2}=x
_{4n+1},u_{4n+3
}=x_{4n+2},u_0=x_{-1},
u_1=x_0,v_0=y_{-1}$  and $v_1=y_0$ we get the result
\bn \ba x_{4n-1}=\frac{x_{-1}}{(1-x_{-1}^2y_0^2)^n},\,  x_{4n}=\frac{x_0(1-y_{-1}x_0)^{2n}}{(1-2y_{-1}x_0)^n},\\
x_{4n+1}=\frac{x_{-1}}{(1+x_{-1}y_0)(1-x_{-1}^2y_0^2)^n}, \,
x_{4n+2}=\frac{(-1)^n(y_{-1}x_0-1)^{2n+1}}{(2y_{-1}x_0-1)^{n+1}},\\
y_{4n-1}=\frac{y_{-1}(1-2y_{-1}x_0)^n}{(1-y_{-1}x_0)^{2n}},\,
y_{4n}=y_0(1-x_{-1}^2y_0^2)^n,\\
y_{4n+1}=\frac{-y_{-1}(1-2y_{-1}x_0)^n}{(1-y_{-1}x_0)^{2n+1}},\,
y_{4n+2}=-y_0(1+x_{-1}y_0)(1-x_{-1}^2y_0^2)^n\ea\en
{obtained by Elsayed \cite{elsa}}. His restriction ($x_{-1},y_{-1},x_0$ and $y_0$ are nonzero real numbers) coincides with our restriction ${x_0y_{-1}\neq 0}$ and ${y_0x_{-1}\neq 0}$ in this case.
\item If { $b=1,a\neq 1$}, then the solutions \eqref{eq35} become
\bn\ba u_{2n}&=u_0 \prod\limits_{r=0}^{n-1}\left\lbrace \frac{a^r [1-a-2u_0v_1]+2u_0v_1}{a^{r+1} [1-a-2u_0v_1]+u_0v_1(1+a)}\right\rbrace \\
v_{2n}&=v_0 \prod\limits_{r=0}^{n-1}\left\lbrace \frac{a^r [1-a-v_0u_1(1+a)]+v_0u_1(1+a)}{a^{r} [1-a-v_0u_1(1+a)]+2v_0u_1} \right\rbrace  \\
u_{2n+1}&=u_1(1-a)\frac{\prod\limits_{r=0}^{n-1}\left\lbrace a^{r} [1-a-v_0u_1(1+a)]+2v_0u_1\right\rbrace }{\prod\limits_{r=0}^{n}\left\lbrace a^r [1-a-v_0u_1(1+a)]+v_0u_1(1+a)\right\rbrace }\\
v_{2n+1}&=v_1(1-a)\frac{\prod\limits_{r=0}^{n-1}\left\lbrace a^{r+1} [1-a-2u_0v_1]+u_0v_1(1+a)\right\rbrace }{\prod\limits_{r=0}^{n}\left\lbrace a^r [1-a-2u_0v_1]+2u_0v_1\right\rbrace }.\ea \label{eq39}\en
{For a particular value  $a=-1$, we obtain }
\bn\ba v_{2n}=v_0 \prod\limits_{r=0}^{n-1}\left\lbrace \frac{(-1)^r}{(-1)^r+v_0u_1}\right\rbrace ,
u_{2n}=u_0 \prod\limits_{r=0}^{n-1}\left\lbrace \frac{(-1)^r [1-u_0v_1]+u_0v_1}{(-1)^{r+1} [1-u_0v_1]} \right\rbrace , \\
v_{2n+1}=v_1\frac{\prod\limits_{r=0}^{n-1}\left\lbrace (-1)^{r+1} [1-u_0v_1]\right\rbrace }{\prod\limits_{r=0}^{n}\left\lbrace (-1)^r [1-u_0v_1]+u_0v_1\right\rbrace },
u_{2n+1}=u_1\frac{\prod\limits_{r=0}^{n-1}\left\lbrace (-1)^r +v_0u_1\right\rbrace }{\prod\limits_{r=0}^{n}\left\lbrace (-1)^r \right\rbrace }\ea \label{eq40}\en
which can be split into
\bn \ba v_{4n}=\frac{v_0}{(1-v_0^2u_1^2)^n},\qquad  v_{4n+1}=\frac{v_1(1-u_0v_1)^{2n}}{(1-2u_0v_1)^n},\\
v_{4n+2}=\frac{v_0}{(1+v_0u_1)(1-v_0^2u_1^2)^n}, \qquad
v_{4n+3}=\frac{(-1)^n(u_0v_1-1)^{2n+1}}{(2u_0v_1-1)^{n+1}},\\
u_{4n}=\frac{u_0(1-2u_0v_1)^n}{(1-u_0v_1)^{2n}},\qquad
u_{4n+1}=u_1(1-v_0^2u_1^2)^n,\\
u_{4n+2}=\frac{-u_0(1-2u_0v_1)^n}{(1-u_0v_1)^{2n+1}},\qquad u_{4n+3}=-u_1(1+v_0u_1)(1-v_0^2u_1^2)^n.\label{eq41}\ea\en
If we let $v_{4n}=y_{4n-1},v_{4n+1}=y_{4n},v_{4n+2}=y_{4n+1},v_{4n+3
}=y_{4n+2},u_{4n}=x_{4n-1},u_{4n+1}=x_{4n},u_{4n+2}=x
_{4n+1},u_{4n+3
}=x_{4n+2},u_0=x_{-1},
u_1=x_0,v_0=y_{-1}$  and $v_1=y_0$ we get the result
\bn \ba y_{4n-1}=\frac{y_{-1}}{(1-y_{-1}^2x_0^2)^n},\,  y_{4n}=\frac{y_0(1-x_{-1}y_0)^{2n}}{(1-2x_{-1}y_0)^n},\\
y_{4n+1}=\frac{y_{-1}}{(1+y_{-1}x_0)(1-y_{-1}^2x_0^2)^n}, \,
y_{4n+2}=\frac{(-1)^n(x_{-1}y_0-1)^{2n+1}}{(2x_{-1}y_0-1)^{n+1}},\\
x_{4n-1}=\frac{x_{-1}(1-2x_{-1}y_0)^n}{(1-x_{-1}y_0)^{2n}},\,
x_{4n}=x_0(1-y_{-1}^2x_0^2)^n,\\
x_{4n+1}=\frac{-x_{-1}(1-2x_{-1}y_0)^n}{(1-x_{-1}y_0)^{2n+1}},\,
x_{4n+2}=-x_0(1+y_{-1}x_0)(1-y_{-1}^2x_0^2)^n\ea\en
{obtained by Elsayed \cite{elsa}}. His restriction ($x_{-1},y_{-1},x_0$ and $y_0$ are nonzero real numbers) coincides with our restriction ${x_0y_{-1}\neq 0}$ and ${y_0x_{-1}\neq 0}$ in this case.
\end{itemize}
\subsection{The case $ab=1$}
\begin{itemize}
\item If {  $a=b=1$}, from  \eqref{eq29} and using the fact that $s_0=\frac{1}{v_0u_1}$ and $t_0=\frac{1}{u_0v_1}$, we obtain
\begin{align} &T_{2r}=\frac{1+2r u_0v_1}{u_0v_1} , \, S_{2r+1}=\frac{1+(2r+1)u_0v_1}{u_0v_1},\,
S_{2r}=\frac{1+2r v_0u_1}{v_0u_1},\,\nonumber\\&T_{2r+1}=\frac{1+(2r+1)v_0u_1}{v_0u_1}.\end{align}
The solutions of the system \eqref{eq2} in this case are
\bn\ba u_{2n}=u_0 \prod\limits_{r=0}^{n-1}\left\lbrace \frac{1+2r u_0v_1}{1+(2r+1)u_0v_1}\right\rbrace  ,&u_{2n+1}=u_1\frac{\prod\limits_{r=0}^{n-1}\left\lbrace 1+(2r+1)v_0u_1 \right\rbrace }{\prod\limits_{r=0}^{n}\left\lbrace 1+2rv_0u_1\right\rbrace }\\
v_{2n}=v_0 \prod\limits_{r=0}^{n-1}\left\lbrace \frac{1+2r v_0u_1}{1+(2r+1)v_0u_1}\right\rbrace  ,&v_{2n+1}=v_1\frac{\prod\limits_{r=0}^{n-1}\left\lbrace 1+(2r+1)u_0v_1 \right\rbrace }{\prod\limits_{r=0}^{n}\left\lbrace 1+2ru_0v_1\right\rbrace }.\ea \en
If we 'shift back', i.e, we let $u_{2n}=x_{2n-1}, v_{2n}=y_{2n-1},u_{2n+1}=x_{2n},v_{2n+1}=y_{2n},u_0=x_{-1},u_1=x_0,v_0=y_{-1}$  and $v_1=y_0$, we get the result
 \bn\ba x_{2n-1}=x_{-1} \prod\limits_{r=0}^{n-1}\left\lbrace \frac{1+2r x_{-1}y_0}{1+(2r+1)x_{-1}y_0}\right\rbrace  ,&x_{2n}=x_0\frac{\prod\limits_{r=0}^{n-1}\left\lbrace 1+(2r+1)y_{-1}x_0 \right\rbrace }{\prod\limits_{r=0}^{n}\left\lbrace 1+2ry_{-1}x_0\right\rbrace }\\
y_{2n-1}=y_{-1} \prod\limits_{r=0}^{n-1}\left\lbrace \frac{1+2r y_{-1}x_0}{1+(2r+1)y_{-1}x_0}\right\rbrace  ,&y_{2n}=v_1\frac{\prod\limits_{r=0}^{n-1}\left\lbrace 1+(2r+1)x_{-1}y_0 \right\rbrace }{\prod\limits_{r=0}^{n}\left\lbrace 1+2rx_{-1}y_0\right\rbrace }\ea \en
{obtained by Kurbanli  et al. \cite{kur}}. The restrictions made by the authors ($x_{-1},y_{-1},x_0$ and $y_0$ are positive real numbers) are included in our
restriction  $-\frac{1}{y_{-1}x_{0}}\neq \{0,2,...,2n\},\quad-\frac{1}{x_{-1}y_{0}}\neq \{0,2,...,2n\},\quad-\frac{1}{x_{-1}y_{0}}\neq \{0,1,...,2n-1\}$ and $-\frac{1}{y_{-1}x_{0}}\neq \{0,1,...,2n-1\}$ in this case.
\item If { $a=b=-1$}, \eqref{eq29} become
\bn \ba T_{2r}=\frac{1}{u_0v_1} ,S_{2r+1}=\frac{-1+u_0v_1}{u_0v_1},
S_{2r}=\frac{1}{v_0u_1} ,T_{2r+1}=\frac{-1+v_0u_1}{v_0u_1}\ea \en
and  solutions of system \eqref{eq2} in this case are
\bn\ba u_{2n}=&\frac{u_0}{(u_0v_1-1)^n} , v_{2n}=\frac{v_0}{(v_0u_1-1)^n},
u_{2n+1}=u_1(v_0u_1-1)^n ,\\ v_{2n+1}=&v_1(u_0v_1-1)^n.\ea \en
\end{itemize}
\section{Example 2}
Consider the third order system of difference equations \eqref{eq112} \cite{elsaibra}.  Imposing the linearized  symmetry condition \eqref{eq4},
we obtain after a set of long calculations the characteristics
\begin{eqnarray}
Q_1=(C_2(-1)^n-C_1 )x_n\nonumber \text{ and }  Q_2=(C_1 +C_2(-1)^n)y_n.
\nonumber\end{eqnarray}
The corresponding generators of the Lie point symmetry are
\bn
\begin{array}{ll}
X_1=-x_n \p_{x_n}+y_n \p_{y_n} \text{ and }
X_2=(-1)^ny_n \p_{u_n}+(-1)^ny_n \p_{y_n}.
\end{array}\label{eq1113}
\en
 The characteristic equation
$$\frac{\d x_n}{x_n}=\frac{\d y_n}{y_n}=\frac{\d x_{n+1}}{-x_{n+1}}=\frac{\d y_{n+1}}{-y_{n+1}}=\frac{\d x_{n+2}}{x_{n+2}}=\frac{\d y_{n+2}}{y_{n+2}}$$
gives us all the invariants under $X_2$  in \eqref{eq1113}. Two of them are given by
\bny w_n=x_n y_{n+1}\label{eq1114}\\z_n=y_nx_{n+1}. \label{eq1115}\eny
From \eqref{eq1114} and \eqref{eq1115},  we have
\bny w_{n+2}=\frac{z_n}{c+d z_n},\quad z_{n+2}=\frac{w_n}{a +b w_n}. \nonumber\eny
By the transformations
\bn z_n=\frac{1}{T_n}, \quad w_n=\frac{1}{S_n}, \label{eq1116}\en
we obtain the linear system
\bn
\ba  S_{n+2}=c T_n +d ,\quad T_{n+2}=a S_n +b .\label{eq1117} \ea\en
which solutions are given by
\begin{subequations}\label{eq1118}
\bn \ba S_n=& \frac{(1+i^n+(-i)^n+(-1)^n)}{4}\left[ a^{\frac{n}{4}}c^{\frac{n}{4}}S_0 + d \sum\limits_{i=0}^{\frac{n-4}{4}}(ac)^i+bc \sum\limits_{i=0}^{\frac{n-4}{4}}(ac)^i\right]
+\\&\frac{(1+i(i^n)-i(-i)^n-(-1)^n)}{4}\left[ a^{\frac{n-3}{4}}c^{\frac{n+1}{4}}T_1 + d \sum\limits_{i=0}^{\frac{n-3}{4}}(ac)^i+bc \sum\limits_{i=0}^{\frac{n-7}{4}}(ac)^i\right]
+\\&\frac{(1-(i^n)-(-i)^n+(-1)^n)}{4}\left[ a^{\frac{n-2}{4}}c^{\frac{n+2}{4}}T_0 + d \sum\limits_{i=0}^{\frac{n-2}{4}}(ac)^i+bc \sum\limits_{i=0}^{\frac{n-6}{4}}(ac)^i\right]
+\\&\frac{(1-i(i^n)+i(-i)^n-(-1)^n)}{4}\left[ a^{\frac{n-1}{4}}c^{\frac{n-1}{4}}S_1 + d \sum\limits_{i=0}^{\frac{n-5}{4}}(ac)^i+bc \sum\limits_{i=0}^{\frac{n-5}{4}}(ac)^i\right]\\
\ea\en\bn\ba
T_n=& \frac{(1+i^n+(-i)^n+(-1)^n)}{4}\left[ a^{\frac{n}{4}}c^{\frac{n}{4}}T_0 + b \sum\limits_{i=0}^{\frac{n-4}{4}}(ac)^i+ad \sum\limits_{i=0}^{\frac{n-4}{4}}(ac)^i\right]
+\\&
\frac{(1+i(i^n)-i(-i)^n-(-1)^n)}{4}\left[ a^{\frac{n+1}{4}}c^{\frac{n-3}{4}}S_1 + b \sum\limits_{i=0}^{\frac{n-3}{4}}(ac)^i+ad \sum\limits_{i=0}^{\frac{n-7}{4}}(ac)^i\right]
+\\&\frac{(1-(i^n)-(-i)^n+(-1)^n)}{4}\left[ a^{\frac{n
+2}{4}}c^{\frac{n-2}{4}}S_0 + b \sum\limits_{i=0}^{\frac{n-2}{4}}(ac)^i+ad \sum\limits_{i=0}^{\frac{n-6}{4}}(ac)^i\right]
+\\&
\frac{(1-i(i^n)+i(-i)^n-(-1)^n)}{4}\left[ a^{\frac{n-1}{4}}c^{\frac{n-1}{4}}T_1 + b \sum\limits_{i=0}^{\frac{n-5}{4}}(ac)^i+ad \sum\limits_{i=0}^{\frac{n-5}{4}}(ac)^i\right].
\\
\ea\en
\end{subequations}
The latter solutions can be split as follows:
\bn \ba S_{4n}=a^nc^n S_0 + d \sum\limits_{i=0}^{{n-1}}(ac)^i+bc \sum\limits_{i=0}^{{n-1}}(ac)^i,\;
T_{4n}=a^nc^n T_0 + b \sum\limits_{i=0}^{{n-1}}(ac)^i\\+ad \sum\limits_{i=0}^{{n-1}}(ac)^i;\,
S_{4n+1}=a^nc^n S_1 + d \sum\limits_{i=0}^{{n-1}}(ac)^i+bc \sum\limits_{i=0}^{{n-1}}(ac)^i,\,
T_{4n+1}=a^nc^n T_1 +\\ b \sum\limits_{i=0}^{{n-1}}(ac)^i+ad \sum\limits_{i=0}^{{n-1}}(ac)^i;\,S_{4n+2}=a^nc^{n+1} T_0 + d \sum\limits_{i=0}^{{n}}(ac)^i+bc \sum\limits_{i=0}^{{n-1}}(ac)^i,\,\\
T_{4n+2}=a^{n+1}c^n S_0 + b \sum\limits_{i=0}^{{n}}(ac)^i+ad \sum\limits_{i=0}^{{n-1}}(ac)^i;\,
S_{4n+3}=a^nc^{n+1} T_1 + d \sum\limits_{i=0}^{{n}}(ac)^i\\+bc \sum\limits_{i=0}^{{n-1}}(ac)^i,\,
T_{4n+3}=a^{n+1}c^n S_1 + b \sum\limits_{i=0}^{{n}}(ac)^i+ad \sum\limits_{i=0}^{{n-1}}(ac)^i.\label{eq1119} \ea\en
From \eqref{eq1114} and \eqref{eq1115}, we have
\bn \ba x_{n+1}=\frac{z_n}{y_n},\quad y_{n+1}=\frac{w_n}{x_n}. \ea \en
By invoking \eqref{eq1116}, the latter system becomes
\bn \ba x_{n+1}=\frac{1}{T_n y_n},\quad y_{n+1}=\frac{1}{S_nx_n}, \label{eq1121}\ea \en
where the functions $S_n$ and $T_n$ are given in \eqref{eq1118}.\\
{\bf Note}: The order of the system \eqref{eq112} has been reduced by two.\\
The  solution of \eqref{eq1121} is given
\bn \ba x_n=\frac{x_0}{2}\left[ \prod\limits_{r=0}^{\frac{n-2}{2}}\frac{S_{2r}}{T_{2r+1}}+(-1)^n\prod\limits_{r=0}^{\frac{n-2}{2}}\frac{S_{2r}}{T_{2r+1}} \right] +\frac{1}{2y_0}\left[\frac{\prod\limits_{r=0}^{\frac{n-3}{2}}S_{2r+1}}{\prod\limits_{r=0}^{\frac{n-1}{2}}T_{2r}}-(-1)^n \frac{\prod\limits_{r=0}^{\frac{n-3}{2}}S_{2r+1}}{\prod\limits_{r=0}^{\frac{n-1}{2}}S_{2r}}
 \right] \\
 y_n=\frac{y_0}{2}\left[ \prod\limits_{r=0}^{\frac{n-2}{2}}\frac{T_{2r}}{S_{2r+1}}+(-1)^n\prod\limits_{r=0}^{\frac{n-2}{2}}\frac{T_{2r}}{S_{2r+1}} \right] +\frac{1}{2x_0}\left[\frac{\prod\limits_{r=0}^{\frac{n-3}{2}}T_{2r+1}}{\prod\limits_{r=0}^{\frac{n-1}{2}}S_{2r}}-(-1)^n \frac{\prod\limits_{r=0}^{\frac{n-3}{2}}T_{2r+1}}{\prod\limits_{r=0}^{\frac{n-1}{2}}S_{2r}}
 \right]. \label{eq1132} \ea \en
The solutions \eqref{eq1132} can be split into
\bn \ba x_{2n}=x_0 \prod\limits_{r=0}^{n-1}\frac{S_{2r}}{T_{2r+1}};y_{2n}=y_0 \prod\limits_{r=0}^{n-1}\frac{T_{2r}}{S_{2r+1}};x_{2n+1}=\frac{\prod\limits_{r=0}^{n-1}
S_{2r+1}}{y_0\prod\limits_{r=0}^{n}T_{2r}};
y_{2n+1}=\frac{\prod\limits_{r=0}^{n-1}T_{2r+1}}{x_0\prod\limits_{r=0}^{n}S_{2r}}. \label{eq1133}\ea \en
Invoking \eqref{eq1118}, the latter solutions become
\bn\ba x_{2n}=x_0 \prod\limits_{r=0}^{n-1}\\\left\lbrace\frac{\alpha\left( a^{\frac{r}{2}}c^{\frac{r}{2}}S_0+d\sum\limits_{i=0}^{\frac{r-2}{2}}(ac)^i+bc\sum\limits_{i=0}^{\frac{r-2}{2}}(ac)^i\right) +\beta\left( a^{\frac{r-1}{2}}c^{\frac{r+1}{2}}T_0+d\sum\limits_{i=0}^{\frac{r-1}{2}}(ac)^i +bc\sum\limits_{i=0}^{\frac{r-3}{2}}(ac)^i\right) }
{\gamma\left( a^{\frac{r+1}{2}}c^{\frac{r-1}{2}}S_1+b\sum\limits_{i=0}^{\frac{r-1}{2}}(ac)^i+ad\sum\limits_{i=0}^{\frac{r-3}{2}}(ac)^i\right) +\lambda\left( a^{\frac{r}{2}}c^{\frac{r}{2}}T_1+b\sum\limits_{i=0}^{\frac{r-2}{2}}(ac)^i +ad\sum\limits_{i=0}^{\frac{r-2}{2}}(ac)^i\right)} \right\rbrace \\
y_{2n}=y_0\prod\limits_{r=0}^{n-1}\\\left\lbrace \frac{\alpha\left( a^{\frac{r}{2}}c^{\frac{r}{2}}T_0+b\sum\limits_{i=0}^{\frac{r-2}{2}}(ac)^i+ad\sum\limits_{i=0}^{\frac{r-2}{2}}(ac)^i\right) +\beta\left( a^{\frac{r+1}{2}}c^{\frac{r-1}{2}}S_0+b\sum\limits_{i=0}^{\frac{r-1}{2}}(ac)^i +ad\sum\limits_{i=0}^{\frac{r-3}{2}}(ac)^i\right)}{\gamma\left( a^{\frac{r-1}{2}}c^{\frac{r+1}{2}}T_1+d\sum\limits_{i=0}^{\frac{r-1}{2}}(ac)^i+bc\sum\limits_{i=0}^{\frac{r-3}{2}}(ac)^i\right) +\lambda\left( a^{\frac{r}{2}}c^{\frac{r}{2}}S_1+d\sum\limits_{i=0}^{\frac{r-2}{2}}(ac)^i +bc\sum\limits_{i=0}^{\frac{r-2}{2}}(ac)^i\right)} \right\rbrace \\
x_{2n+1}=\\\frac{\prod\limits_{r=0}^{n-1}\left\lbrace \gamma\left( a^{\frac{r-1}{2}}c^{\frac{r+1}{2}}T_1+d\sum\limits_{i=0}^{\frac{r-1}{2}}(ac)^i+bc\sum\limits_{i=0}^{\frac{r-3}{2}}(ac)^i\right) +\lambda\left( a^{\frac{r}{2}}c^{\frac{r}{2}}S_1+d\sum\limits_{i=0}^{\frac{r-2}{2}}(ac)^i +bc\sum\limits_{i=0}^{\frac{r-2}{2}}(ac)^i\right)\right\rbrace }{y_0\prod\limits_{r=0}^{n}\left\lbrace\alpha\left( a^{\frac{r}{2}}c^{\frac{r}{2}}T_0+b\sum\limits_{i=0}^{\frac{r-2}{2}}(ac)^i+ad\sum\limits_{i=0}^{\frac{r-2}{2}}(ac)^i\right) +\beta\left( a^{\frac{r+1}{2}}c^{\frac{r-1}{2}}S_0+b\sum\limits_{i=0}^{\frac{r-1}{2}}(ac)^i +ad\sum\limits_{i=0}^{\frac{r-3}{2}}(ac)^i\right)\right\rbrace } \\
y_{2n+1}=\\\frac{\prod\limits_{r=0}^{n-1}\left\lbrace \gamma\left( a^{\frac{r+1}{2}}c^{\frac{r-1}{2}}S_1+b\sum\limits_{i=0}^{\frac{r-1}{2}}(ac)^i+ad\sum\limits_{i=0}^{\frac{r-3}{2}}(ac)^i\right) +\lambda\left( a^{\frac{r}{2}}c^{\frac{r}{2}}T_1+b\sum\limits_{i=0}^{\frac{r-2}{2}}(ac)^i +ad\sum\limits_{i=0}^{\frac{r-2}{2}}(ac)^i\right) \right\rbrace }{x_0\prod\limits_{r=0}^{n}\left\lbrace \alpha\left( a^{\frac{r}{2}}c^{\frac{r}{2}}S_0+d\sum\limits_{i=0}^{\frac{r-2}{2}}(ac)^i+bc\sum\limits_{i=0}^{\frac{r-2}{2}}(ac)^i\right) +\beta\left( a^{\frac{r-1}{2}}c^{\frac{r+1}{2}}T_0+d\sum\limits_{i=0}^{\frac{r-1}{2}}(ac)^i +bc\sum\limits_{i=0}^{\frac{r-3}{2}}(ac)^i\right)  \right\rbrace},\label{eq1124}\ea \en
where $\alpha=(1+(-1)^r)/2,\,\beta=(1-(1)^r)/2,\,\gamma=(1-(-1)^r)/2$ and $\lambda=(1+(1)^r)/2$. The initials conditions are
\bn S_0=\frac{1}{x_0y_1},\quad S_1=\frac{1}{x_1y_2},\quad T_0=\frac{1}{y_0x_1},\quad T_1=\frac{1}{y_1x_2}.\en
We get the solution to the system \eqref{eq112} by splitting  \eqref{eq1124}. We have
\begin{subequations}\label{eq1126}
\bn \ba

 x_{4n}=\frac{(x_2y_2)^n}{y_0^nx_0^{n-1}}
\left\lbrace \frac{\splitfrac{
\prod\limits_{r=0}^{n-1}\big\{ a^{r}c^{r}+dx_0y_1 \sum\limits_{i=0}^{r-1}(ac)^i+bc x_0y_1\sum\limits_{i=0}^{r-1}(ac)^i\big\} \prod\limits_{r=0}^{n-1}\big\{ a^{r}c^{r+1}}{+dy_0x_1 \sum\limits_{i=0}^{r}(ac)^i+bc y_0x_1\sum\limits_{i=0}^{r-1}(ac)^i\big\}}
}
{\splitfrac{\prod\limits_{r=0}^{n-1}\big\{ a^{r+1}c^{r}+bx_1y_2 \sum\limits_{i=0}^{r}(ac)^i+ad x_1y_2\sum\limits_{i=0}^{r-1}(ac)^i\big\} \prod\limits_{r=0}^{n-1}\big\{ a^{r}c^{r}}{+by_1x_2 \sum\limits_{i=0}^{r-1}(ac)^i+ad y_1x_2\sum\limits_{i=0}^{r-1}(ac)^i\big\}}}\right\rbrace
\ea \en
\bn \ba
x_{4n+2}=\frac{x_2^{n+1}y_2^n}{y_0^nx_0^{n}}\left\lbrace\frac{\splitfrac{\prod\limits_{r=0}^{n}\big\{ a^{r}c^{r}+dx_0y_1 \sum\limits_{i=0}^{r-1}(ac)^i+bc x_0y_1\sum\limits_{i=0}^{r-1}(ac)^i\big\} \prod\limits_{r=0}^{n-1}\big\{ a^{r}c^{r+1}+}{dy_0x_1 \sum\limits_{i=0}^{r}(ac)^i+bc y_0x_1\sum\limits_{i=0}^{r-1}(ac)^i\big\} } }{\splitfrac{\prod\limits_{r=0}^{n-1}\big\{ a^{r+1}c^{r}+bx_1y_2 \sum\limits_{i=0}^{r}(ac)^i+ad x_1y_2\sum\limits_{i=0}^{r-1}(ac)^i\big\} \prod\limits_{r=0}^{n}\big\{ a^{r}c^{r}+}{by_1x_2 \sum\limits_{i=0}^{r-1}(ac)^i+ad y_1x_2\sum\limits_{i=0}^{r-1}(ac)^i\big\} }}\right\rbrace
\ea \en
\bn \ba
y_{4n}=\frac{(x_2y_2)^n}{x_0^ny_0^{n-1}}\left\lbrace\frac{\splitfrac{\prod\limits_{r=0}^{n-1}\big\{ a^{r}c^{r}+by_0x_1 \sum\limits_{i=0}^{r-1}(ac)^i+ad y_0x_1\sum\limits_{i=0}^{r-1}(ac)^i\big\} \prod\limits_{r=0}^{n-1}\big\{ a^{r+1}c^{r}}{+bx_0y_1 \sum\limits_{i=0}^{r}(ac)^i+ad x_0y_1\sum\limits_{i=0}^{r-1}(ac)^i\big\}  }}{\splitfrac{\prod\limits_{r=0}^{n-1}\big\{ a^{r}c^{r+1}+dy_1x_2 \sum\limits_{i=0}^{r}(ac)^i+bc y_1x_2\sum\limits_{i=0}^{r-1}(ac)^i\big\} \prod\limits_{r=0}^{n-1}\big\{ a^{r}c^{r}}{+dx_1y_2 \sum\limits_{i=0}^{r-1}(ac)^i+bc x_1y_2\sum\limits_{i=0}^{r-1}(ac)^i\big\} }}\right\rbrace
\ea \en
\bn \ba
y_{4n+2}=\frac{x_2^ny_2^{n+1}}{x_0^ny_0^{n-1}}\left\lbrace\frac{\splitfrac{\prod\limits_{r=0}^{n}\big\{ a^{r}c^{r}+by_0x_1 \sum\limits_{i=0}^{r-1}(ac)^i+ad y_0x_1\sum\limits_{i=0}^{r-1}(ac)^i\big\} \prod\limits_{r=0}^{n-1}\big\{ a^{r+1}c^{r}}{+bx_0y_1 \sum\limits_{i=0}^{r}(ac)^i+ad x_0y_1\sum\limits_{i=0}^{r-1}(ac)^i\big\} } }{\splitfrac{\prod\limits_{r=0}^{n-1}\big\{ a^{r}c^{r+1}+dy_1x_2 \sum\limits_{i=0}^{r}(ac)^i+bc y_1x_2\sum\limits_{i=0}^{r-1}(ac)^i\big\} \prod\limits_{r=0}^{n}\big\{ a^{r}c^{r}}{+dx_1y_2 \sum\limits_{i=0}^{r-1}(ac)^i+bc x_1y_2\sum\limits_{i=0}^{r-1}(ac)^i\big\}} }\right\rbrace
 \ea \en
\bn \ba
 x_{4n+1}=\frac{x_1 (y_0x_0)^n}{(x_2y_2)^n} \left\lbrace\frac{\splitfrac{\prod\limits_{r=0}^{n-1}\big\{ a^{r}c^{r+1}+dy_1x_2 \sum\limits_{i=0}^{r}(ac)^i+bc y_1x_2\sum\limits_{i=0}^{r-1}(ac)^i\big\} \prod\limits_{r=0}^{n-1}\big\{ a^{r}c^{r}}{+dx_1y_2 \sum\limits_{i=0}^{r-1}(ac)^i+bc x_1y_2\sum\limits_{i=0}^{r-1}(ac)^i\big\}}}{\splitfrac{\prod\limits_{r=0}^{n}\big\{ a^{r}c^{r}+by_0x_1 \sum\limits_{i=0}^{r-1}(ac)^i+ad y_0x_1\sum\limits_{i=0}^{r-1}(ac)^i\big\}\prod\limits_{r=0}^{n-1}\big\{ a^{r+1}c^{r}}{+bx_0y_1 \sum\limits_{i=0}^{r}(ac)^i+ad x_0y_1\sum\limits_{i=0}^{r-1}(ac)^i\big\}}}\right\rbrace
 \ea \en
\bn \ba
 x_{4n+3}=\frac{y_1 (y_0x_0)^n}{x_2^ny_2^{n+1}} \left\lbrace\frac{\splitfrac{\prod\limits_{r=0}^{n-1}\big\{ a^{r}c^{r+1}+dy_1x_2 \sum\limits_{i=0}^{r}(ac)^i+bc y_1x_2\sum\limits_{i=0}^{r-1}(ac)^i\big\} \prod\limits_{r=0}^{n}\big\{ a^{r}c^{r}+}{dx_1y_2 \sum\limits_{i=0}^{r-1}(ac)^i+bc x_1y_2\sum\limits_{i=0}^{r-1}(ac)^i\big\}}}{\splitfrac{\prod\limits_{r=0}^{n}\big\{ a^{r}c^{r}+by_0x_1 \sum\limits_{i=0}^{r-1}(ac)^i+ad y_0x_1\sum\limits_{i=0}^{r-1}(ac)^i\big\} \prod\limits_{r=0}^{n}\big\{ a^{r+1}c^{r}+}{bx_0y_1 \sum\limits_{i=0}^{r}(ac)^i+ad x_0y_1\sum\limits_{i=0}^{r-1}(ac)^i\big\}}}\right\rbrace
\ea \en
\bn \ba
y_{4n+1}=\frac{y_1 (x_0y_0)^n}{(x_2y_2)^n}\left\lbrace\frac{\splitfrac{\prod\limits_{r=0}^{n-1}\big\{ a^{r+1}c^{r}+bx_1y_2 \sum\limits_{i=0}^{r}(ac)^i+ad x_1y_2\sum\limits_{i=0}^{r-1}(ac)^i\big\} \prod\limits_{r=0}^{n-1}\big\{ a^{r}c^{r}+}{by_1x_2 \sum\limits_{i=0}^{r-1}(ac)^i+ad y_1x_2\sum\limits_{i=0}^{r-1}(ac)^i\big\}}}{\splitfrac{\prod\limits_{r=0}^{n}\big\{ a^{r}c^{r}+dx_0y_1 \sum\limits_{i=0}^{r-1}(ac)^i+bc x_0y_1\sum\limits_{i=0}^{r-1}(ac)^i\big\} \prod\limits_{r=0}^{n-1}\big\{ a^{r}c^{r+1}}{+dy_0x_1 \sum\limits_{i=0}^{r}(ac)^i+bc y_0x_1\sum\limits_{i=0}^{r-1}(ac)^i\big\} }}\right\rbrace
\ea \en
\bn \ba
y_{4n+3}=\frac{x_1 x_0^ny_0^{n+1}}{x_2^{n+1}y_2^n}\left\lbrace\frac{\splitfrac{\prod\limits_{r=0}^{n-1}\big\{ a^{r+1}c^{r}+bx_1y_2 \sum\limits_{i=0}^{r}(ac)^i+ad x_1y_2\sum\limits_{i=0}^{r-1}(ac)^i\big\} \prod\limits_{r=0}^{n}\big\{ a^{r}c^{r}}{+by_1x_2 \sum\limits_{i=0}^{r-1}(ac)^i+ad y_1x_2\sum\limits_{i=0}^{r-1}(ac)^i\big\}}}{\splitfrac{\prod\limits_{r=0}^{n}\big\{ a^{r}c^{r}+dx_0y_1 \sum\limits_{i=0}^{r-1}(ac)^i+bc x_0y_1\sum\limits_{i=0}^{r-1}(ac)^i\big\} \prod\limits_{r=0}^{n}\big\{ a^{r}c^{r+1}}{+dy_0x_1 \sum\limits_{i=0}^{r}(ac)^i+bc y_0x_1\sum\limits_{i=0}^{r-1}(ac)^i\big\} }}\right\rbrace
\ea\en
\end{subequations}
with
\begin{align} & \frac{1}{x_0y_1} \neq \frac{d+bc}{-a^rc^r} \sum\limits_{i=0}^{{r-1}}(ac)^i;
 \frac{1}{y_0x_1} \neq \frac{b+ad}{-a^rc^r} \sum\limits_{i=0}^{{r-1}}(ac)^i;\nonumber\\
 &\frac{1}{x_1y_2} \neq\frac{d+bc}{-a^rc^r}  \sum\limits_{i=0}^{{r-1}}(ac)^i;
\frac{1}{y_1x_2} \neq \frac{b+ad}{-a^rc^r } \sum\limits_{i=0}^{{r-1}}(ac)^i ;\nonumber
\\ & 1 \neq \frac{y_0x_1}{-a^rc^{r+1}}\left[ d \sum\limits_{i=0}^{{r}}(ac)^i+bc \sum\limits_{i=0}^{{r-1}}(ac)^i\right] ;
 1 \neq\frac{x_0y_1}{-a^{r+1}c^r}  \left[b \sum\limits_{i=0}^{{r}}(ac)^i+ad \sum\limits_{i=0}^{{r-1}}(ac)^i\right], \nonumber\\&
1 \neq \frac{y_1x_2}{-a^rc^{r+1}}\left[ d \sum\limits_{i=0}^{{r}}(ac)^i+bc \sum\limits_{i=0}^{{r-1}}(ac)^i\right] ;
1 \neq\frac{x_1y_2}{-a^{r+1}c^r}\left[  b \sum\limits_{i=0}^{{r}}(ac)^i+ad \sum\limits_{i=0}^{{r-1}}(ac)^i\right], \label{eq1128}\end{align}
$r\leq n$ and \begin{align} x_0y_1\neq 0,\quad x_1y_2\neq 0,\quad y_0x_1\neq 0 , \quad y_1x_2\neq 0 . \label{eq1127}\end{align}
\subsection{Special cases}
In this section, we investigate some special cases which  lead us to all results found in \cite{elsaibra}.\\ 
$\bullet $ If {$ac\neq 1$}, then the solutions \eqref{eq1126} become
\begin{subequations}\label{eq1129}
\bn\ba x_{4n}=\frac{(x_2y_2)^n}{x_0^{n-1}y_0^n}\\\frac{\prod\limits_{r=0}^{n-1}\left\lbrace a^{r}c^{r}[1-ac-x_0y_1(d+bc)]+x_0y_1(d+bc) \right\rbrace \prod\limits_{r=0}^{n-1}\left\lbrace a^{r}c^{r}[c-ac^2-y_0x_1(acd+bc)]+y_0x_1(d+bc) \right\rbrace }{\prod\limits_{r=0}^{n-1}\left\lbrace a^{r}c^{r}[a-a^2c-x_1y_2(abc+ad)]+x_1y_2(b+ad) \right\rbrace \prod\limits_{r=0}^{n-1}\left\lbrace a^{r}c^{r}[1-ac-y_1x_2(b+ad)]+y_1x_2(b+ad) \right\rbrace } \\
\ea \en
\bn\ba
x_{4n+2}=\frac{x_2^{n+1}y_2^n}{x_0^{n}y_0^n}\\\frac{\prod\limits_{r=0}^{n}\left\lbrace a^{r}c^{r}[1-ac-x_0y_1(d+bc)]+x_0y_1(d+bc) \right\rbrace \prod\limits_{r=0}^{n-1}\left\lbrace a^{r}c^{r}[c-ac^2-y_0x_1(acd+bc)]+y_0x_1(d+bc) \right\rbrace }{\prod\limits_{r=0}^{n-1}\left\lbrace a^{r}c^{r}[a-a^2c-x_1y_2(abc+ad)]+x_1y_2(b+ad) \right\rbrace \prod\limits_{r=0}^{n}\left\lbrace a^{r}c^{r}[1-ac-y_1x_2(b+ad)]+y_1x_2(b+ad) \right\rbrace }
\ea\en
\bn \ba
y_{4n}=\frac{(x_2y_2)^n}{x_0^{n}y_0^{n-1}}\\\frac{\prod\limits_{r=0}^{n-1}\left\lbrace a^{r}c^{r}[1-ac-y_0x_1(b+ad)]+y_0x_1(b+ad) \right\rbrace \prod\limits_{r=0}^{n-1}\left\lbrace a^{r}c^{r}[a-a^2c-x_0y_1(abc+ad)]+x_0y_1(b+ad) \right\rbrace }{\prod\limits_{r=0}^{n-1}\left\lbrace a^{r}c^{r}[c-ac^2-y_1x_2(acd+bc)]+y_1x_2(d+bc) \right\rbrace \prod\limits_{r=0}^{n-1}\left\lbrace a^{r}c^{r}[1-ac-x_1y_2(d+bc)]+x_1y_2(d+bc) \right\rbrace }
\ea\en
\bn \ba
y_{4n+2}=\frac{x_2^ny_2^{n+1}}{x_0^{n}y_0^{n}}\\\frac{\prod\limits_{r=0}^{n}\left\lbrace a^{r}c^{r}[1-ac-y_0x_1(b+ad)]+y_0x_1(b+ad) \right\rbrace \prod\limits_{r=0}^{n-1}\left\lbrace a^{r}c^{r}[a-a^2c-x_0y_1(abc+ad)]+x_0y_1(b+ad) \right\rbrace }{\prod\limits_{r=0}^{n-1}\left\lbrace a^{r}c^{r}[c-ac^2-y_1x_2(acd+bc)]+y_1x_2(d+bc) \right\rbrace \prod\limits_{r=0}^{n}\left\lbrace a^{r}c^{r}[1-ac-x_1y_2(d+bc)]+x_1y_2(d+bc) \right\rbrace }
\ea\en
\bn \ba
x_{4n+1} =\frac{x_1x_0^{n}y_0^{n}(1-ac)}{x_2^{n}y_2^{n}} \\ \frac{\prod\limits_{r=0}^{n-1}\left\lbrace a^{r}c^{r}[c-ac^2-y_1x_2(acd+bc)]+y_1x_2(d+bc) \right\rbrace \prod\limits_{r=0}^{n-1}\left\lbrace a^{r}c^{r}[1-ac-x_1y_2(d+bc)]+x_1y_2(d+bc) \right\rbrace}{\prod\limits_{r=0}^{n}\left\lbrace a^{r}c^{r}[1-ac-y_0x_1(b+ad)]+y_0x_1(b+ad) \right\rbrace \prod\limits_{r=0}^{n-1}\left\lbrace a^{r}c^{r}[a-a^2c-x_0y_1(abc+ad)]+x_0y_1(b+ad) \right\rbrace}
\ea\en
\bn \ba

x_{4n+3} =\frac{y_1x_0^{n+1}y_0^{n}(1-ac)}{x_2^{n}y_2^{n+1}} \\ \frac{\prod\limits_{r=0}^{n-1}\left\lbrace a^{r}c^{r}[c-ac^2-y_1x_2(acd+bc)]+y_1x_2(d+bc) \right\rbrace \prod\limits_{r=0}^{n}\left\lbrace a^{r}c^{r}[1-ac-x_1y_2(d+bc)]+x_1y_2(d+bc) \right\rbrace}{\prod\limits_{r=0}^{n}\left\lbrace a^{r}c^{r}[1-ac-y_0x_1(b+ad)]+y_0x_1(b+ad) \right\rbrace \prod\limits_{r=0}^{n}\left\lbrace a^{r}c^{r}[a-a^2c-x_0y_1(abc+ad)]+x_0y_1(b+ad) \right\rbrace}\\

y_{4n+1}=\frac{y_1(x_0y_0)^n(1-ac)}{(x_2y_2)^n}\\\frac{\prod\limits_{r=0}^{n-1}\left\lbrace a^{r}c^{r}[a-a^2c-x_1y_2(abc+ad)]+x_1y_2(b+ad) \right\rbrace \prod\limits_{r=0}^{n-1}\left\lbrace a^{r}c^{r}[1-ac-y_1x_2(b+ad)]+y_1x_2(b+ad) \right\rbrace }{\prod\limits_{r=0}^{n}\left\lbrace a^{r}c^{r}[1-ac-x_0y_1(d+bc)]+x_0y_1(d+bc) \right\rbrace \prod\limits_{r=0}^{n-1}\left\lbrace a^{r}c^{r}[c-ac^2-y_0x_1(acd+bc)]+y_0x_1(d+bc) \right\rbrace}\\

y_{4n+3}=\frac{x_1x_0^ny_0^{n+1}(1-ac)}{x_2^{n+1}y_2^n}\\\frac{\prod\limits_{r=0}^{n-1}\left\lbrace a^{r}c^{r}[a-a^2c-x_1y_2(abc+ad)]+x_1y_2(b+ad) \right\rbrace \prod\limits_{r=0}^{n}\left\lbrace a^{r}c^{r}[1-ac-y_1x_2(b+ad)]+y_1x_2(b+ad) \right\rbrace }{\prod\limits_{r=0}^{n}\left\lbrace a^{r}c^{r}[1-ac-x_0y_1(d+bc)]+x_0y_1(d+bc) \right\rbrace \prod\limits_{r=0}^{n}\left\lbrace a^{r}c^{r}[c-ac^2-y_0x_1(acd+bc)]+y_0x_1(d+bc) \right\rbrace}.
\ea \en
\end{subequations}
We recover the results found in \cite{elsaibra}
by substituting with some specific values $b$ and $d$. For instance, if we let $b=d=1$, the solutions \eqref{eq1129} become
\begin{subequations}\label{eq1133}
\bn \ba x_{8n}=\frac{(x_2y_2)^{2n}}{x_0^{2n-1}y_0^{2n}}\frac{[1-y_0x_1]^{2n} }{[1+ x_1y_2 ]^{n} [-1+x_1y_2 ]^{n}[-1+2y_1x_2]^{n}} ;&\\

x_{8n+1} =\frac{x_1x_0^{2n}y_0^{2n}}{x_2^{2n}y_2^{2n}}\frac{ [1-y_1x_2]^{2n}}{[1+x_0y_1]^{n}  [-1+2y_0x_1]^{n}[-1+x_0y_1]^{n}}

\ea\en
\bn \ba
x_{8n+2}=\frac{x_2^{2n+1}y_2^{2n}}{x_0^{2n}y_0^{2n}}\frac{[1-y_0x_1]^{2n} }{[1+ x_1y_2 ]^{n} [-1+x_1y_2 ]^{n}[-1+2y_1x_2]^{n}}  ;&\\

x_{8n+3} =\frac{ y_1x_0^{2n+1}y_0^{2n}}{x_2^{2n}y_2^{2n+1}} \frac{ [1-y_1x_2]^{2n}}{[1+x_0y_1]^{n+1}  [-1+2y_0x_1]^{n}[-1+x_0y_1]^{n}}

\ea\en
\bn \ba
x_{8n+4}=\frac{(x_2y_2)^{2n+1}[-1+y_0x_1]^{2n+1} }{x_0^{2n}y_0^{2n+1}[1+ x_1y_2 ]^{n+1} [-1+x_1y_2 ]^{n}[-1+2y_1x_2]^{n}}  ;&\\
x_{8n+5} = \frac{x_1(x_0y_0)^{2n+1} [-1+y_1x_2]^{2n+1}}{(x_2y_2)^{2n+1}[1+x_0y_1]^{n+1}  [-1+2y_0x_1]^{n+1}[-1+x_0y_1]^{n}}

\ea\en
\bn \ba

x_{8n+6}=\frac{x_2^{2n+2}y_2^{2n+1} [1-y_0x_1]^{2n+1} }{x_0^{2n+1}y_0^{2n+1}[1+ x_1y_2 ]^{n+1} [x_1y_2 -1]^{n}[2y_1x_2-1]^{n+1}};&\\
x_{8n+7} = \frac{y_1x_0^{2n+2}y_0^{2n+1} [1-y_1x_2]^{2n+1}}{x_2^{2n+1}y_2^{2n+2}[1+x_0y_1]^{n+1}  [2y_0x_1-1]^{n+1}[-1+x_0y_1]^{n+1}}
\\
\ea\en
\bn \ba
y_{8n}=\frac{(x_2y_2)^{2n}}{x_0^{2n}y_0^{2n-1}}\frac{[1+x_0y_1]^{n}  [-1+2y_0x_1]^{n}[-1+x_0y_1]^{n}}{ [1-y_1x_2]^{2n} };&\\
y_{8n+1}=\frac{y_1(x_0y_0)^{2n}}{(x_2y_2)^{2n}}\frac{[1+ x_1y_2 ]^{n} [-1+x_1y_2 ]^{n}[-1+2y_1x_2]^{n}}{[1-y_0x_1]^{2n} }\\
y_{8n+2}=\frac{x_2^{2n}y_2^{2n+1}}{x_0^{2n}y_0^{2n}}\frac{[1+x_0y_1]^{n}  [-1+2y_0x_1]^{n}[-1+x_0y_1]^{n}}{ [1-y_1x_2]^{2n} };
\ea\en
\bn \ba
y_{8n+3}=\frac{x_1x_0^{2n}y_0^{2n+1}}{x_2^{2n+1}y_2^{2n}}\frac{[1+ x_1y_2 ]^{n} [-1+x_1y_2 ]^{n}[-1+2y_1x_2]^{n}}{[1-y_0x_1]^{2n} }\\
%
y_{8n+4}=\frac{(x_2y_2)^{2n+1}}{x_0^{2n+1}y_0^{2n}}\frac{[1+x_0y_1]^{n+1}  [-1+2y_0x_1]^{n}[-1+x_0y_1]^{n}}{[-1+y_1x_2]^{2n+1}  };
\ea\en
\bn \ba
y_{8n+5}=\frac{y_1(x_0y_0)^{2n+1}[1+ x_1y_2 ]^{n+1} [-1+x_1y_2 ]^{n}[-1+2y_1x_]^{n}}{(x_2y_2)^{2n+1}[1-y_0x_1]^{2n+1} }\\
y_{8n+6}=\frac{x_2^{2n+1}y_2^{2n+2}[1+x_0y_1]^{n+1}  [2y_0x_1-1]^{n+1}[x_0y_1-1]^{n}}{x_0^{2n+1}y_0^{2n+1} [1-y_1x_2]^{2n+1} };&\\
y_{8n+7}=\frac{x_1x_0^{2n+1}y_0^{2n+2}[1+ x_1y_2 ]^{n+1} [x_1y_2-1 ]^{n}[2y_1x_2-1]^{n+1}}{x_2^{2n+2}y_2^{2n+1}[1-y_0x_1]^{2n+2} }.
\ea\en
\end{subequations}
By \lq shifting back \rq twice,
 we get exactly the results found in \cite{elsaibra}. The restrictions made by the authors ($x_{-2},x_{-1},x_{0},y_{-2},y_{-1}$ and $y_{0}$ non zero real number with $x_{-1}y_{0},x_{-2}y_{-1}\neq \pm1$, and $x_0y_{-1},x_{-1}y_{-2}\neq 1,\frac{1}{2}$) are included in our restrictions given \eqref{eq1127} and \eqref{eq1128}.\par \noindent
 $\bullet$ If {$ac=1$}, then the solutions \eqref{eq1126} become
 \begin{subequations}\label{eq1135}
\bn \ba x_{4n}=\frac{(x_2y_2)^n\prod\limits_{r=0}^{n-1}\left\lbrace 1+x_0y_1 r(d+bc)\right\rbrace \prod\limits_{r=0}^{n-1}\left\lbrace c +dy_0x_1 (r+1)+bc y_0x_1r\right\rbrace  }{y_0^nx_0^{n-1}\prod\limits_{r=0}^{n-1}\left\lbrace a +bx_1y_2 (r+1)+ad x_1y_2r \right\rbrace \prod\limits_{r=0}^{n-1}\left\lbrace 1+y_1x_2 r(b+ad)\right\rbrace };\\

x_{4n+2}=\frac{x_2^{n+1}y_2^n\prod\limits_{r=0}^{n}\left\lbrace 1+x_0y_1 r(d+bc)\right\rbrace \prod\limits_{r=0}^{n-1}\left\lbrace c+y_0x_1 (dr+d+bcr)+\right\rbrace  }{y_0^nx_0^{n}\prod\limits_{r=0}^{n-1}\left\lbrace a+x_1y_2 (br+b+adr)\right\rbrace \prod\limits_{r=0}^{n}\left\lbrace 1+y_1x_2 r(b+ad)\right\rbrace }
\ea\en
\bn \ba
y_{4n}=\frac{(x_2y_2)^n\prod\limits_{r=0}^{n-1}\left\lbrace 1+by_0x_1 r+ad y_0x_1r\right\rbrace \prod\limits_{r=0}^{n-1}\left\lbrace a+x_0y_1 (br+b+adr)\right\rbrace  }{x_0^ny_0^{n-1}\prod\limits_{r=0}^{n-1}\left\lbrace c+dy_1x_2 (r+1)+bc y_1x_2r\right\rbrace \prod\limits_{r=0}^{n-1}\left\lbrace 1+dx_1y_2 r+bc x_1y_2r\right\rbrace };\\

y_{4n+2}=\frac{x_2^ny_2^{n+1}\prod\limits_{r=0}^{n}\left\lbrace 1+y_0x_1 r(b+ad)\right\rbrace \prod\limits_{r=0}^{n-1}\left\lbrace a+x_0y_1 (br+b+adr)\right\rbrace  }{x_0^ny_0^{n-1}\prod\limits_{r=0}^{n-1}\left\lbrace c+y_1x_2 (dr+d+bcr)\right\rbrace \prod\limits_{r=0}^{n}\left\lbrace 1+x_1y_2 r(d+bc)\right\rbrace }
\ea\en
\bn \ba

 x_{4n+1}= \frac{x_1 (y_0x_0)^n\prod\limits_{r=0}^{n-1}\left\lbrace c+y_1x_2 (dr+d+bcr)\right\rbrace \prod\limits_{r=0}^{n-1}\left\lbrace 1+x_1y_2 r(d+bc)\right\rbrace}{(x_2y_2)^n\prod\limits_{r=0}^{n}\left\lbrace 1+y_0x_1 r(b+ad)\right\rbrace \prod\limits_{r=0}^{n-1}\left\lbrace a+x_0y_1 (br+b+adr)\right\rbrace};\\

 x_{4n+3}=\frac{y_1 (y_0x_0)^n}{x_2^ny_2^{n+1}} \frac{\prod\limits_{r=0}^{n-1}\left\lbrace c +y_1x_2 (dr+d+bcr)\right\rbrace \prod\limits_{r=0}^{n}\left\lbrace 1+x_1y_2 r(d+bc) \right\rbrace}{\prod\limits_{r=0}^{n}\left\lbrace 1+y_0x_1 r(b+ad)\right\rbrace \prod\limits_{r=0}^{n}\left\lbrace a+x_0y_1 (br+b+adr)\right\rbrace}
 \ea\en
\bn \ba

y_{4n+1}=\frac{y_1 (x_0y_0)^n\prod\limits_{r=0}^{n-1}\left\lbrace a+x_1y_2 (br+b+adr)\right\rbrace \prod\limits_{r=0}^{n-1}\left\lbrace 1+y_1x_2 r(b+ad)\right\rbrace}{(x_2y_2)^n\prod\limits_{r=0}^{n}\left\lbrace 1+x_0y_1 r(d+bc)\right\rbrace \prod\limits_{r=0}^{n-1}\left\lbrace c+y_0x_1 (dr+d+bcr)\right\rbrace };\\

y_{4n+3}=\frac{x_1 x_0^ny_0^{n+1}}{x_2^{n+1}y_2^n}\frac{\prod\limits_{r=0}^{n-1}\left\lbrace a+x_1y_2 (br+b+adr)\right\rbrace \prod\limits_{r=0}^{n}\left\lbrace 1+y_1x_2 r(b+ad)  \right\rbrace}{\prod\limits_{r=0}^{n}\left\lbrace 1+x_0y_1 r(d+bc)\right\rbrace \prod\limits_{r=0}^{n}\left\lbrace c+y_0x_1 (dr+d+bcr)\right\rbrace }
\label{eq1134}\ea\en
We recover results found in \cite{elsaibra} by letting $a,b,c$ and $d$ take specific values keeping $ac=1$. For instance, if $a=b=c=d=1$, we get the solutions
\bn \ba x_{4n}=\frac{(x_2y_2)^n\prod\limits_{r=0}^{n-1}\left\lbrace 1+x_0y_1 2r\right\rbrace \prod\limits_{r=0}^{n-1}\left\lbrace 1 +y_0x_1 (2r+1)\right\rbrace  }{y_0^nx_0^{n-1}\prod\limits_{r=0}^{n-1}\left\lbrace 1 +x_1y_2 (2r+1) \right\rbrace \prod\limits_{r=0}^{n-1}\left\lbrace 1+y_1x_2 2r\right\rbrace }\ea\en
\bn \ba

x_{4n+2}=\frac{x_2^{n+1}y_2^n\prod\limits_{r=0}^{n}\left\lbrace 1+x_0y_1 2r\right\rbrace \prod\limits_{r=0}^{n-1}\left\lbrace 1+y_0x_1 (2r+1)\right\rbrace  }{y_0^nx_0^{n}\prod\limits_{r=0}^{n-1}\left\lbrace 1+x_1y_2 (2r+1)\right\rbrace \prod\limits_{r=0}^{n}\left\lbrace 1+y_1x_2 2r\right\rbrace }
\ea\en
\bn \ba

y_{4n}=\frac{(x_2y_2)^n\prod\limits_{r=0}^{n-1}\left\lbrace 1+y_0x_1 2r\right\rbrace \prod\limits_{r=0}^{n-1}\left\lbrace 1+x_0y_1 (2r+1)\right\rbrace  }{x_0^ny_0^{n-1}\prod\limits_{r=0}^{n-1}\left\lbrace 1+y_1x_2 (2r+1)\right\rbrace \prod\limits_{r=0}^{n-1}\left\lbrace 1+x_1y_2 2r\right\rbrace }\ea\en
\bn \ba

y_{4n+2}=\frac{x_2^ny_2^{n+1}\prod\limits_{r=0}^{n}\left\lbrace 1+y_0x_1 2r\right\rbrace \prod\limits_{r=0}^{n-1}\left\lbrace 1+x_0y_1 (2r+1)\right\rbrace  }{x_0^ny_0^{n-1}\prod\limits_{r=0}^{n-1}\left\lbrace 1+y_1x_2 (2r+1)\right\rbrace \prod\limits_{r=0}^{n}\left\lbrace 1+x_1y_2 2r\right\rbrace }\\

 x_{4n+1}= \frac{x_1 (y_0x_0)^n\prod\limits_{r=0}^{n-1}\left\lbrace 1+y_1x_2 (2r+1)\right\rbrace \prod\limits_{r=0}^{n-1}\left\lbrace 1+x_1y_2 2r\right\rbrace}{(x_2y_2)^n\prod\limits_{r=0}^{n}\left\lbrace 1+y_0x_1 2r\right\rbrace \prod\limits_{r=0}^{n-1}\left\lbrace a+x_0y_1 (2r+1)\right\rbrace}\ea\en
\bn \ba

 x_{4n+3}=\frac{y_1 (y_0x_0)^n}{x_2^ny_2^{n+1}} \frac{\prod\limits_{r=0}^{n-1}\left\lbrace 1 +y_1x_2 (2r+1)\right\rbrace \prod\limits_{r=0}^{n}\left\lbrace 1+x_1y_2 2r \right\rbrace}{\prod\limits_{r=0}^{n}\left\lbrace 1+y_0x_1 2r\right\rbrace \prod\limits_{r=0}^{n}\left\lbrace 1+x_0y_1 (2r+1)\right\rbrace}
 \ea\en
\bn \ba

y_{4n+1}=\frac{y_1 (x_0y_0)^n\prod\limits_{r=0}^{n-1}\left\lbrace 1+x_1y_2 (2r+1)\right\rbrace \prod\limits_{r=0}^{n-1}\left\lbrace 1+y_1x_2 2r\right\rbrace}{(x_2y_2)^n\prod\limits_{r=0}^{n}\left\lbrace 1+x_0y_1 2r\right\rbrace \prod\limits_{r=0}^{n-1}\left\lbrace c+y_0x_1 (2r+1)\right\rbrace }\ea\en
\bn \ba

y_{4n+3}=\frac{x_1 x_0^ny_0^{n+1}}{x_2^{n+1}y_2^n}\frac{\prod\limits_{r=0}^{n-1}\left\lbrace 1+x_1y_2 (2r+1)\right\rbrace \prod\limits_{r=0}^{n}\left\lbrace 1+y_1x_2 2r  \right\rbrace}{\prod\limits_{r=0}^{n}\left\lbrace 1+x_0y_1 2r\right\rbrace \prod\limits_{r=0}^{n}\left\lbrace 1+y_0x_1 (2r+1)\right\rbrace }.
\ea\en
\end{subequations}
By "shifting back" the solutions in \eqref{eq1135},
we get exactly the results found in \cite{elsaibra}. The restrictions made by the authors ($x_{-2},x_{-1},x_{0},y_{-2},y_{-1}$ and $y_{0}$ non zero real number ) are exactly  our restrictions given \eqref{eq1127}.
\section{Conclusion}
We have  used symmetry methods to reduce the order of systems of difference equations  and have found their  exact solutions. Solutions of special cases of the systems investigated  exist in recent literature ( see Elsayed \cite{elsa}, Kurbanli et al. \cite{kur}, \cite{elsaibra}). These authors used induction methods obtain their results. For the sake of clarification, we have shifted back  our solutions to show that their solutions  match the results obtained in this work.

\end{document}